\pdfoutput=1 
\RequirePackage{fix-cm}
\documentclass{amsart}
\usepackage{fixltx2e}
\RequirePackage{eucal}
\RequirePackage{optparams}
\RequirePackage{enumerate}
\RequirePackage{ifmtarg}
\RequirePackage{xifthen}
\RequirePackage{suffix}
\RequirePackage{amsmath}
\RequirePackage{amssymb}
\RequirePackage{amsthm}
\RequirePackage{amsfonts}
\RequirePackage{onlyamsmath}   
\usepackage{undertilde}
\usepackage{hyperxmp}
\usepackage{hyperref}
\hypersetup{		pdfborder=0 0 0,%
								colorlinks=true,%
								linkbordercolor={0 0 0},
								hyperfootnotes=true,%
								pdftex=true, %
								implicit=true,%
								final=true,%
								letterpaper=true, %
								bookmarks=true,                       
                bookmarksnumbered=true,               
                bookmarksopen=true,							%
								pdfpagemode=UseOutlines,                
	              pdfstartview=Fit,                
	              pdfpagelayout=OneColumn,            
	              breaklinks=true,                 
	              pageanchor=true,                 
	              plainpages=false,                
	              pdfpagelabels=true,                   
	              backref=page,%
	              pagebackref=true,
	              hyperindex=true,                      
								pdfdisplaydoctitle=true,%
								pdfauthor={Peter M. Gerdes},%
pdfkeywords={REA, w-REA, Turing reducible, computability},%
pdfsubject={MSC2010: Primary 03D25, 03D30; Secondary (03D60)},%
pdftitle={A w-REA Set Forming A Minimal Pair With 0'}}

\newtheorem{theorem}{Theorem}[section]
\newtheorem{lemma}[theorem]{Lemma}
\newtheorem{proposition}[theorem]{Proposition}

\theoremstyle{definition}
\newtheorem{definition}[theorem]{Definition}

\theoremstyle{remark}

\numberwithin{equation}{section}

\makeatletter

\newcommand*{\abs}[1]{\lvert#1\rvert}
\newcommand*{\defword}[1]{\textbf{#1}}
\newcommand*{\ensuretext}[1]{\ensuremath{\text{#1}} }
\DeclareMathOperator{\dom}{dom}
\providecommand*{\nin}{\mathrel{\not\in}}
\newcommand*{\conv}[1][]{\mathpunct{\downarrow}_{#1}}
\newcommand*{\diverge}{\mathpunct{\uparrow}}
\newcommand*{\powset}[1]{\mathcal{P}\left(#1\right)}
\newcommand*{\restr}[1]{\mathpunct{\restriction_{#1}}}

\providecommand*{\pair}[2]{\mathopen{\langle} #1, #2 \mathclose{\rangle}}
\newcommand*{\setcol}[2]{{#1}^{[#2]}}
\newcommand*{\lh}[2][]{\abs{#2}_{#1}}
\newcommand*{\incompat}{\mathrel{\mid}}
\newcommand*{\compat}{\mathrel{\nmid}}
\newcommand*{\bstrs}{2^{<\omega}}
\newcommand*{\wstrs}{\omega^{<\omega}}
\newcommand*{\deltaZeroN}[1]{\Delta^{0}_{#1}}
\newcommand*{\deltaZeroTwo}{\deltaZeroN{2}}
\WithSuffix\def\deltaZeroN[#1]#2{\deltaZeroN{#2}\!\left(#1\right)}
\newcommand*{\set}[2]{\ifthenelse{\isempty{#2}}{\left \{  #1 \right \}}{\left \{ #1 \middle | #2\right \}} }
\newcommand*{\jump}[1]{#1'}
\newcommand*{\eset}{\emptyset}
\newcommand*{\Tdeg}[1]{\utilde{#1}}
\newcommand*{\Tzerosym}{0}
\newcommand*{\Tzero}{\Tdeg{\Tzerosym}}
\newcommand*{\zeroj}{\jump{\Tzerosym}}
\newcommand*{\zerojj}{\jump{\jump{\Tzerosym}}}
\long\def\REset@[#1][#2]#3{W^{#1}_{\ifthenelse{\isempty{#2}}{#3}{#3,#2}}}
\newcommand{\REset}{\optparams{\REset@}{[][]}}
\newcommand*{\recfnlSYM}{\Phi}
\newcommand*{\recfnl}[4][]{\recfnlSYM_{#2\ifthenelse{\isempty{#1}}{}{,#1} }%
	\ifthenelse{\isempty{#4}}%
	{\ifthenelse{\isempty{#3}}%
		{}
		{(#3)}
	}
	{\ifthenelse{\isempty{#3}}{(\eset;#4)}{(#3;#4)}%
}}
\newcommand{\requirement}[2]{\ensuremath{\mathcal{#1}_{#2}}}
\newcommand{\req}[2]{\requirement{#1}{#2}}
\newcommand*{\re}{\ensuretext{c.e.}}
\newcommand*{\ce}{\ensuretext{c.e.}}
\newcommand*{\wck}{\ensuremath{\omega^{ck}_1}}
\newcommand*{\union}{\mathbin{\cup}}
\newcommand*{\isect}{\mathbin{\cap}}

\newenvironment{reqs}{%
	\renewcommand{\req}[2]{\item[\llap{\ensuremath{\displaystyle \requirement{##1}{##2}}\hphantom{X}:}]}%
	\vspace{.5em}%
	\begin{description}%
	\addtolength{\itemindent}{4em}
}{\end{description}\vspace{.5em}}

\newcommand*{\nTequiv}{\ensuremath{\ncong_T}}

\newcommand*{\TSYM}{\mathbf{T}}

\newcommand*{\Tleq}{\leq_{\TSYM}}

\newcommand*{\Tgeq}{\geq_{\TSYM}}

\newcommand*{\Tless}{<_{\TSYM}}
\newcommand*{\nTleq}{\nleq_{\TSYM}}

\newcommand*{\Tmeet}{\mathbin{\wedge_{\TSYM}}}
\newcommand*{\Tjoin}{\mathbin{\vee_{\TSYM}}}
\newcommand*{\alphaREA}[1][\alpha]{\ensuremath{ #1 \text{-REA}} }
\newcommand*{\REA}{\ensuremath{\text{REA}}}
\WithSuffix\def\REA[#1]{\ifthenelse{\isempty{#1}}{\alphaREA[\alpha]}{\alphaREA[{#1}]}}
\WithSuffix\def\exists[#1]{\left(\exists\, #1\right)\!}
\WithSuffix\def\forall[#1]{\left(\forall\, #1\right)\!}
\newcommand*{\str}[1]{\mathopen{\langle}#1\mathclose{\rangle}}
\newcommand*{\concat}{\mathbin{\widehat{}}}
\WithSuffix\def\concat[#1]{\concat\str{#1}}
\newcommand*{\setcmp}[1]{\overline{#1}}
\newcommand*{\sigmaZeroN}[1]{\Sigma^0_{#1}}
\newcommand*{\sigmaZeroOne}{\sigmaZeroN{1}}
\newcommand*{\sigmazi}{\sigmaZeroOne}
\newcommand*{\Tplus}{\mathbin{\oplus}}
\WithSuffix\def\Tplus[#1][#2]{\mathop{\bigoplus}^{#2}_{#1}}
\newcommand*{\kleeneO}[1][]{\mathcal{O}^{#1}}
\newcommand*{\kleenel}{<_{\kleeneO}}
\newcommand*{\kleenePlus}{\mathbin{+_{\kleeneO}}}
\newcommand*{\cequiv}{\mathrel{\backsimeq}}
\makeatother

\newcommand{\laxiom}[3]{ \mathopen{\langle}#1: \ifthenelse{\isempty{#2}}{\eset}{#2} \rightarrow #3 \mathclose{\rangle} }
\providecommand*{\subfun}{\prec}
\providecommand*{\supfun}{\succ}

\newcommand{\axset}{\mathcal{A}}

\begin{document}

	\title{A \( \omega \)-REA Set Forming A Minimal Pair With \( 0' \)}
	\author{Peter M. Gerdes}
	\address{Department of Mathematics\\
University of Notre Dame du Lac\\
Notre Dame, Indiana 46556}
	\urladdr{\href{http://invariant.org/}{http://invariant.org/}}
	\email{\href{mailto:gerdes@invariant.org}{gerdes@invariant.org}}
	\thanks{Partially Supported by NSF EMSW21-RTG-0739007 and EMSW21-RTG-0838506}
	\date{in Fall 2010 for consideration for publication in the \href{http://www.aslonline.org/journals-journal.html}{Journal for Symbolic Logic}}

\begin{abstract}
	It is easy to see that no \( n \)-REA set can form a (non-trivial) minimal pair with \( \zeroj \) and only slightly more difficult to observe that no \( \omega \)-REA set can form a (non-trivial) minimal pair with \( \zerojj \).  Shore has asked whether this can be improved to show that no \( \omega \)-REA set forms a (non-trivial) minimal pair with \( \zeroj \).  We show that no such improvement is possible by constructing a set \( C \) with \( \Tzero \Tless C \Tleq \zerojj \) forming a minimal pair with \( \zeroj \).  We then show that no \( \alpha \)-REA set can form a (non-trivial) minimal pair with \( \zerojj \).
\end{abstract}

\maketitle

%

\section{Introduction}

\subsection{Notation}

The notation we use in this proof is largely standard.  We use \( \sigma, \tau, \delta \) to denote partial functions from \( \omega \) to \( \set{0,1}{} \) and write \( \sigma \subfun \tau \) to denote that the function \( \tau \) extends \( \sigma \).  We identify sets with their characteristic functions so that \( \sigma \subfun X \) has the expected meaning for \( X \in \powset{\omega} \).   We denote \( x \in \dom \sigma \) ( \( x \nin \dom \sigma \) ) by \( \sigma(x)\conv \) ( \( \sigma(x)\diverge \) ) and say \( \sigma \) is incompatible (compatible) with \( \tau \), denoted \( \sigma \incompat \tau \) (\( \sigma \compat \tau \)), if there is some (no) \( x \) with \( \sigma(x)\conv \neq \tau(x)\conv  \).  We let \( \alpha, \beta, \gamma \) range over elements in \( \wstrs \), write \( \alpha^{-} \) as shorthand for \( \alpha\restr{\lh{\alpha}-1} \) and denote the concatenation of \( \alpha \) with \( \beta \) by \( \alpha\concat\beta \).  We denote the length of a \( \alpha \) by \( \lh{\alpha} \) and the and extend this notation to partial functions by setting \( \lh{\delta} \) = \( 1 + \max \dom \delta \).  Capital roman letters range over subsets of \( \omega \) which we identify with their characteristic function.  

\( \recfnl{e}{Z}{x} \) denotes the \( e \)-th \( \set{0,1}{} \) valued partial computable functional applied to oracle \( Z \) on the input \( x \).  We adopt the convention that if \( \recfnl{e}{Z}{x} \) converges in  \( s \) steps, written \( \recfnl{e}{Z}{x}\conv[s] \), then \( \recfnl{e}{Z\restr{s}}{x}=\recfnl{e}{Z}{x} \).  \( \REset[Z]{e} \) is the \( e \)-th set \re in \( Z \) and  \( \REset[Z][s]{e} \) is it's stage \( s \) approximation.  We use \( \pair{x}{y} \) to denote the integer code of the pair \( (x, y)  \).  Capital roman letters range over \( \powset{\omega} \) and we write  \( \Tdeg{C} \) for the Turing degree of \( C \), \( \setcmp{C} \) for the compliment of \( C \), \( \jump{C} \) for the jump of \( C \) and use \( \Tleq, T\equiv, \Tmeet, \Tjoin \) to denote Turing reducibility, equivalence, meet and join respectively.  

We follow the standard practice of identifying \( \setcol{X}{n} \), the \( n \)-th column of \( X \), with \( \set{y }{\pair{n}{y} \in X} \) and  \( \setcol{X}{\leq n} \) for \( \set{ \pair{m}{y} \in X }{m \leq n} \).  We extend this notation to partial functions by letting \( \setcol{\sigma}{\leq l} \) represent the restriction of \( \sigma \) to \( \setcol{\omega}{\leq l} \).  We also stipulate that \( \setcol{X}{<0} = \eset \).

\subsection{Overview}

In \cite{pseudojump-operators-II} Jockusch and Shore introduce the \( \REA[\alpha] \), for \( \alpha < \wck \), sets as the sets produced by effectively iterating the construction of a relatively \re set \( \alpha \) many times.  Since we will restrict our attention here to \( \alpha = \omega \) we will use the equivalent (up to Turing degree) definition.

\begin{definition}\label{def:w-rea}
	\( C \subseteq \omega \) is \( \REA[\omega] \)  iff there is a computable function \( f \) such that

	\begin{equation*}
		\setcol{C}{n}= \REset[{\setcol{C}{<n}}]{f(n)}
	\end{equation*}
\end{definition}

Professor Shore has observed that if \( C \nTleq \Tzero \) is \( \REA[\omega] \) then \( \exists B \Tleq C \) with \( \Tzero \Tless B \Tless \zerojj \) but asked (private communication) if this would still hold if \( \zerojj \) was replaced with \( \zeroj \).  In this paper we answer this question in the negative by proving the following theorem.

\begin{theorem}\label{thm:no-delta2-bounding}
	There is an \( \REA[\omega] \) set \( C \nTequiv \Tzero \) such that \( \zeroj \Tmeet \Tdeg{C} = \Tzero \)
\end{theorem}

\subsection{Failure at \( \zerojj \)}

Before we embark on this construction it is instructive to see why this claim fails for \( \zerojj \).

\begin{proposition}\label{prop:delta3-bounding}
	If \( C \) is \( \REA[\omega] \) and \( C \nTequiv \Tzero \) then \( C \Tmeet \zerojj \neq \Tzero \).
\end{proposition}
\begin{proof}

Assume \( C \) fails the lemma.  We first argue that for every \( n \) \( \setcol{C}{\leq n} \) must be computable.  Since \(  \setcol{C}{n+1} \) is \( \deltaZeroTwo \) in \(\setcol{C}{\leq n} \) if \( n \) is greatest with  \( \setcol{C}{\leq n} \Tleq \Tzero \) we must have \( \Tzero \Tless \setcol{C}{n+1} \Tleq \zeroj \Tleq \zerojj \).  Hence if \( C \Tmeet \zerojj = \Tzero \) then \( \setcol{C}{\leq n} \Tleq \Tzero \) for all \( n \in \omega \).

So suppose \( \setcol{C}{\leq n} \Tleq \Tzero \) for all \( n \in \omega \).  We now argue that \( \zerojj \) can compute \( C \).  Note that given a \ce index for a computable set \( R \) \( \zerojj \) can recover a \ce index for \( \setcmp{R} \) and from an index for \( \setcol{C}{\leq n} \) and \( \setcmp{\setcol{C}{\leq n}} \) one can uniformly recover a \ce index from \( \setcol{C}{\leq n+1} \).  Thus if \( \setcol{C}{\leq n} \) is computable for all \( n \in \omega \) by induction \( \zerojj \) can recover \( i_n, j_n \) with \( \setcol{C}{ \leq n} = \REset{i_n} = \setcmp{\REset{j_n}} \).    Clearly these indexes allow \( \zerojj \) to compute \( C \).

\end{proof}

The lesson to be drawn from this proof is that any \( C \) satisfying theorem \ref{thm:no-delta2-bounding} must be the join of a countable collection of computable sets.  Thus the non-computability of \( C \) must result from the non-uniformity of this join.  The difficulty in building \( C \) is therefore how to encode enough about \( \zeroj \) in \( \setcol{C}{<n} \) so \( \setcol{C}{n} \) can successfully diagonalize against the \( \deltaZeroN{2} \) sets while making sure \( \setcol{C}{<n} \) only encodes a finite amount of non-computable information.

\section{Machinery}

\subsection{Building \( \omega \)-REA Sets}

Evidently if we are to build \( C \nin \deltaZeroN{2} \) we will have to somehow have to uniformly specify an \re procedure to build \(  \setcol{C}{n} \) from \( \setcol{C}{<n} \) while dealing with the fact that our approximation to \( \setcol{C}{<n} \) will never settle on the correct value.  Rather than trying to explicitly give such a procedure upfront we will instead enumerate rules called axioms committing us to enumerate certain elements into \( \setcol{C}{n} \) when certain conditions are met by \( \setcol{C}{<n} \).    

\begin{definition}\label{def:re-axiom}
	A \defword{axiom} is a triple  \( \laxiom{l}{\sigma}{y} \) where \( l \in \omega \),  \( \sigma \) is a function from a finite subset of \( \setcol{\omega}{< l} \) to \( \set{0,1}{} \) and \( y \in \setcol{\omega}{\geq l} \).
\end{definition}

In our construction we will think of the axiom \( \laxiom{l}{\sigma}{y} \) as the commitment to place \( y \) in \(  C \) if \( \sigma \subfun C \).  The parameter \( l \) serves only to ensure that attempts to enumerate elements in the \( n \)-th column of \( C \) are only allowed to consult the first \( n-1 \) columns of \( C \) thus avoiding any circularity.  
The utility of this definition is made clear by the following lemma.

\begin{lemma}\label{lem:build-wrea}
	If \( \axset \) is an \re set of axioms then the set \( C \) defined by
	\begin{equation*}
		y \in C \iff \exists[l \in \omega]\exists[\sigma \subfun C ]\Bigl[ \laxiom{l}{\sigma}{y} \in \axset \Bigr]
	\end{equation*}
	is \( \REA[\omega] \)
\end{lemma}
\begin{proof}
	Note that
	\begin{equation*}
		\pair{n}{x} \in C \iff  \exists[l\leq n]\exists[\smash{\sigma \subfun \setcol{C}{< l}}]\Bigl[ \laxiom{l}{\sigma}{y} \in \axset  \Bigr]
	\end{equation*}
	Thus \( \setcol{C}{n} \) only depends on \( \setcol{C}{<n}  \) so \( C \) is well defined.  Furthermore the above equation explicitly defines \( \setcol{C}{n} \) from \( \setcol{C}{<n} \)  and \(  n  \) via a (uniformly) \( \sigmazi \) formula.  Thus by an application of the s-m-n theorem \cite{general-recursive-functions-of-natural-numbers} there is a computable function \( f \) satisfying definition \ref{def:w-rea}.
\end{proof}

Our construction will proceed by building a \re set \( \axset \) of axioms which will yield an \( \REA[\omega] \) set via the preceding lemma.  To make proper use of this machinery we introduce two more definitions.  We first try and capture the notion that some axiom \( \laxiom{l}{\sigma'}{y} \) only has an effect if \( \sigma \subfun C \).

\begin{definition}\label{def:dependent-on}
	The axiom \( \laxiom{l}{\sigma'}{y} \) \defword{depends} on \( \sigma \) if \( \sigma \subfun \sigma' \).  We say the axiom \( \laxiom{l}{\sigma}{y} \) is enumerated dependent on \( \delta \) to mean we enumerate \( \laxiom{l}{ \sigma \union \setcol{\delta}{< l} }{y} \) into \( \axset \).
\end{definition}

We will also speak of an axiom depending on \( C(n)=0 \) to mean it depends on the \( \sigma \) defined by  \( \sigma(n)= 0 \).  During our construction we will frequently want to satisfy some requirement on the assumption that a guess about how \( C \) behaves on some finite number of columns and a finite initial segment is true.  We therefore introduce a notion of how the axioms would affect \( C \) if such a guess were correct. 

\begin{definition}\label{def:yields-over}
	Given any set \( C_{\alpha} \subseteq \setcol{\omega}{< l_\alpha} \) and a partial function \( \delta_{\alpha}  \)  satisfying \( \setcol{\delta_\alpha}{<l} \subfun \setcol{C_\alpha}{<l_\alpha}  \) (understood as a guess at an initial segment of \( C \)) we say that a set of axioms \( \axset \) yields \( C \) over \( C_{\alpha}, \delta_{\alpha} \) if
	\begin{equation*}
		\begin{aligned}
		\pair{n}{x} \in C  \iff &  \biggl( n < l_\alpha \land \pair{n}{x} \in C_\alpha \biggr) \lor \\
														  & \biggl( \pair{n}{x} \in \delta_\alpha \biggr) \lor \\
														& \biggl( \pair{n}{x}  \nin \dom \delta_\alpha \land \exists[l \leq n]\exists[\sigma \subfun C ]\Bigl[ \laxiom{l}{\sigma}{\pair{n}{x}} \in \axset \Bigr] \biggr)
	\end{aligned}
	\end{equation*}
\end{definition}

In other words \( \axset \) yields \( X \) over \( C_{\alpha}, \delta_{\alpha} \) if we take \( C_\alpha, \delta_\alpha \) to be the first \( l_\alpha \) columns of \( C \) and \( \delta_\alpha \subfun C \) regardless of what the axioms say and then build the rest of \( C \) using the construction from lemma \ref{lem:build-wrea}.

\section{Requirements \& Modules}

We fix a computable array \( V_{e,s} \) of finite sets via the limit lemma \cite{limit-lemma} such that every \(  \deltaZeroN{2} \) set is of the form \( V_e = \lim_{s\to\infty} V_{e,s} \) and build \( C \) to meet the following requirements.

\begin{reqs}
	\req{R}{e,i} \( \recfnl{i}{C}{} \neq V_e \) or \( V_e \Tleq \Tzero \) whenever \( V_e \) defined.
	\req{N}{e} \( \REset{e} \neq \setcmp{C}  \)
\end{reqs}

We reserve columns \( 3\pair{e}{i} \), \( 3\pair{e}{i}+1 \) for \( \req{R}{e,i} \) and the column \( 3e+2 \) for meeting \( \req{N}{e} \) and grant each requirement the right to modify a finite initial segment of later columns but not earlier columns.  Each column of \( C \) will be either finite or co-finite thereby making \( \setcol{C}{<n} \) computable as our observation required.

As \( C \) can't be computable in \( \zeroj \) during the construction later requirements won't know, even in the limit, how the earlier requirements are satisfied.  To deal with this we perform our construction along a tree assigning to each \( \alpha \in \wstrs \) in the tree a module \( P_\alpha \) tasked with handling a particular requirement on the assumption that \( \alpha \) correctly encodes how the higher priority requirements are met.  In particular we assign requirements to modules as follows.

\begin{equation}\label{handles}
	P_\alpha \text{ handles } \begin{cases}
																	\req{R}{e,i} & \text{ if } \lh{\alpha} = 2\pair{e}{i} \\
																	\req{N}{e}  & \text{ if } \lh{\alpha} = 2e + 1
															\end{cases}
\end{equation}

Given \( \alpha \in \wstrs \) we define \( l_\alpha \) to be the first column reserved for the requirement handled by \( P_\alpha \).  Note that \( l_\alpha \) only depends on \( \lh{\alpha} \) so all modules tasked with meeting a given requirement share columns.  We will associate to each \( \alpha \) set \( C_{\alpha} \subset \setcol{\omega}{<l_\alpha} \) intended as a guess at \( \setcol{C}{< l_\alpha} \) and a partial function \( \delta_{\alpha} \in \bstrs \) representing a guess at the finite part of \( C \) used by prior requirements.  These two guesses will always be compatible, i.e. \( \setcol{\delta_\alpha}{< l_\alpha} \subfun \setcol{C}{< l_\alpha} \).  Implicitly \( \delta_{\alpha} \) will function as a restraint as well since \( P_{\alpha} \) won't attempt to change \( C(x) \) if \( x \in \dom \delta_\alpha \).  We regard \( C_{\alpha}, \delta_{\alpha} \) as a description of the ultimate effect of \( P_{\alpha^{-}} \) on \( C \).

\subsection{Action Along The Tree}

As explained above the module \( P_\alpha \) will act to meet it's requirement using the information encoded in \( \alpha \) about how earlier requirements were met.  At each stage we will have some guess at how the various requirements are met and that guess will control which modules are then executed at that stage, i.e., only those modules that appear to have correct guesses execute. Those familiar with \( \Pi^0_2 \) tree constructions may be assured that the tree executes the modules in the standard fashion and skip ahead to the next section while those desiring more details can read on. 

More formally we will define a function \( f \), the true path, such that if \( \alpha^{-} \subset f \) then \( f(\lh{\alpha}) \) indicates how \( P_\alpha \) satisfies it's associated requirement.  In an abuse of notation we will write \( f(\alpha) = n \) to indicate that if \( f \supseteq \alpha \) then \( f(\lh{\alpha})=n \).  At any stage \( s \) we will have some approximation \( f_{s} \in \wstrs \) to the true path with \( f = \liminf_{s} f_{s}  \).  We will execute a single module \( P_\alpha \)  satisfying \( f_{s} \supseteq \alpha \) at every stage \( s \) and leave it to \( P_\alpha \) to set the value of \( f_s(\alpha) \) at such stages.  We ensure that if \( f_s \supset \alpha \) occurs infinitely often then \( P_\alpha \) is executed infinitely often as well by starting out at the root node and executing in increasing order the modules at each node \( \alpha \subseteq f_s \) with \( \lh{\alpha} \leq l  \) before starting over at the root and working out to nodes of length \( l + 1 \). 


\section{The Construction}

A full description of the construction will consist of giving the behavior of each module \( P_\alpha \) the approximation to it's outcome \( f_{s}(\alpha) \) and the properties \( C_{\beta}, \delta_{\beta} \) for each \( \beta=\alpha\concat[f_{s}(\alpha)] \).  We will always define \( C_{\beta}, \delta_{\beta} \) at the first stage \( f_s \supseteq \beta \) guaranteeing they are always defined when needed.  Note that when describing the various modules we will say the \( P_\alpha \) stage \( s \) to refer to the \( s \)-th time the module \( P_\alpha \) is executed.  We will also adopt the shorthand \( \alpha^{+} \) for \( \alpha\concat[f(\alpha)] \) whenever \( \alpha \subset f \).

\subsection{Basic Approach}

Before describing the full construction it's useful to informally sketch how each requirement is to be met.   The action of the module \( P_\alpha \) implementing the strategy \( \req{N}{e} \) can be thought of as implementing a straightforward finite injury argument as follows.  \( P_\alpha \) will wait for a chance to enumerate some element from \( W_e \) into \( \setcol{C}{l_\alpha} \) doing nothing until such an element is found.  If no such element is found then both \( W_e \) and \( C \) fail to cover some element in the column \( l_\alpha \).  On the other hand if such an element is found \( \req{N}{e} \) will enumerate that element into \( \setcol{C}{l_\alpha}  \) and reset all weaker priority requirements.  This reset is accomplished simply by permanently changing \( f_{s}(\alpha) \) from the \( 0 \) it had been up till now to \( 1 \) thereby abandoning all previously visited modules \( P_\beta, \beta \supsetneq \alpha \).

The interesting case occurs when \( P_\alpha \) implements \( \req{R}{e,i} \).  Here our strategy will be to lay dormant (unactivated) as long as the action of weaker requirements never leads us to change our mind about (our approximation to) \( \recfnl{i}{C}{} \), i.e., yields only compatible computations.  If we remain in this situation we will argue that \( \recfnl{i}{C}{} \) is computable.  If we do see a change in \( \recfnl{i}{C}{x} \) for some \( x \) we will activate \( \req{R}{e,i} \) and work to alternate between the two computations to diagonalize against \( V_{e}(x) \).  Later we will show that if we ever change our mind about \( \recfnl{i}{C}{x} \) then \( P_\alpha \) has the means to roll back the intervening axioms and recover the previous value of \( \recfnl{i}{C}{x} \) by enumerating some controlling element into \( \setcol{C}{l_\alpha+1} \).  \( P_\alpha \) can now act to ensure that \( \recfnl{i}{C}{x} \) always disagrees with \( V_{e,s}(x) \) by taking said element in and out of \( \setcol{C}{l_\alpha+1} \).  To ensure that \( P_\alpha \) can later change take it back out each time \( P_\alpha \) enumerates the controlling element into \( \setcol{C}{l_\alpha+1} \) it does so dependent on some large number being absent from \( \setcol{C}{l_\alpha} \).  By latter adding this number to \( \setcol{C}{l_\alpha} \), \( P_\alpha \) can effectively cancel it's previous commitment and keep \( \recfnl{i}{C}{x} \neq V_{e,s}(x) \).  

Provided \( V_{e,s}(x) \) eventually settles down this provides no problem.  Each time \( V_{e,s}(x) \) flip-flops we simply set \( f_{s}(\alpha) \) to the next unused value which has the effect of resetting all the subsequence requirements.  However, we must accommodate the possibility this limit fails to exist and somehow prevent those \( P_\beta ,  \beta \supset \alpha \) that assume the limit exists from interfering with those that assume it doesn't.  The key point here is to ensure that a particular flag element will be in \( \setcol{C}{l_\alpha+1} \) iff \( \lim_{s\to\infty} V_{e,s}(x) \) exists.  This allows the modules \( P_\beta, \beta \supset \alpha \) guessing the limit doesn't exist to predicate all their actions on the absence of this element and vice versa ensuring noninterference.  The effect of this is to ensure that if \( \alpha \subset f \) once \( \alpha \) appears on \( f_s \) then no other requirements modify the region of \( C \) used by \( P_\alpha \).

\subsection{Global Constraints}

To ensure the \( P_\alpha \) modules interact appropriately we need to impose two minor additional constraints on the construction.

\begin{enumerate}[(I)]\label{add-const}
	\item If \( P_\alpha \) enumerates axiom \( \pi \) then \( \pi \) is enumerated dependent on \( \delta_\alpha \).\label{add-const:delta_alpha}
	\item If \( P_\alpha \) wants to enumerate axiom \( \pi \) and \( P_\beta, \beta \subsetneq \alpha \) is an unactivated \( \req{R}{e,i} \) module then \( \pi \) is enumerated dependent on the partial function sending \( \pair{l_\alpha +1}{m} \) to \( 0 \) with \( m \) larger than anything mentioned so far in the construction .\label{add-const:reverse}
\end{enumerate}

This first constraint will ensure that if the guess \( \delta_\alpha \) at an initial segment of \( C \) is wrong then the axioms enumerated by \( P_\alpha \) have no effect on the construction.  In particular it will guarantee that if \( P_\beta \) implements \( \req{R}{e,i} \) the modules \( \alpha \supset \beta \) assuming that \( \req{R}{e,i} \) has only finitely many eventful stages and those assuming it has infinitely many such stages don't interfere.

The second constraint will allow \( P_\beta \) implementing \( \req{R}{e,i} \)  to `roll back' axioms enumerated by \( \alpha \supset \beta \) while \( P_\beta \) to recover an earlier computation of \( \recfnl{i}{C}{} \).  This will ensure that even if we first see one value for \( \recfnl{i}{C}{x} \) and then \( \recfnl{i}{C}{x} \) appears to diverge for many stages before converging to an alternate value we will still be able to return to the first value and use it to diagonalize against \( \recfnl{i}{C}{} =V_e \).

We now give the detailed actions of the various modules with the understanding that they be modified in the obvious way to comply with these two constraints. 

\subsection{The Basic \( \req{R}{e,i} \) Strategy }\label{subsec:R_ei-strategy}

Suppose \( P_{\alpha} \) is assigned to handle \( \req{R}{e,i} \).  We wait until we observe a \( P_\alpha \) stage \( s \) (i.e. \( s \)-th time \( P_\alpha \) is executed), integers \( x,y \), strings \( Y_0, Y_1 \in \bstrs \) such that \( C \) would extend \( Y_1 \) if nothing is done but \( Y_0 \) if \( y \nin Y_0 \) is added to \( \setcol{C}{l_\alpha +1 } \) and \( Y_0, Y_1 \) yield incompatible computations.  More formally 

\begin{align}\label{activate-cond}
	\begin{split}
	& \setcol{Y_0}{\leq l_\alpha} = \setcol{Y_1}{\leq l_\alpha} \land \delta_\alpha \subfun Y_0, Y_1 \\
  & \recfnl{i}{Y_0}{x}\conv[y] \neq \recfnl{i}{Y_1}{x}\conv \\
	& \lh{Y_0} < y \\
	& y \in \setcol{\omega}{l_\alpha +1 } \\
  & \axset_s \text{ yields an extension of } Y_1 \text{ over }  C_\alpha, \delta_\alpha \\
  & \axset_s \union \set{ \laxiom{l_\alpha+1}{\eset}{y} }{} \text{ yields an extension of } Y_0 \text{ over } C_\alpha, \delta_\alpha
\end{split}
\end{align}

When this occurs we say that \( P_\alpha \) is \defword{activated} at stage \( s \).  We will later show that these conditions are equivalent to the informal notion of changing our mind about the value of \( \recfnl{i}{C}{x} \).   

If \( P_\alpha \) is activated at \( P_\alpha \) stage \( s \) select \( q \in \setcol{\omega}{l_\alpha +1} \) larger than any number mentioned so far to serve as our `flag' by satisfying \( q \in C  \) iff \( P_\alpha \) only acts finitely many times.  Also pick \( k_{s} \) larger than any number mentioned so far with the intent of enumerating \( \pair{l_\alpha}{k_{t}} \) into \( C \) to cancel any decision to put \( y \) into \( C \) at \( P_\alpha \) stage \( s \).   Now for any \( t \geq s+1 \) let \( j_t \in \set{0,1}{} \) such that  \( \recfnl{i}{Y_{j_t}}{x} \neq V_{e,s}(x)  \).  Say that a \( P_\alpha \) stage \( t > s \) is \defword{eventful} for \( P_\alpha \) if \( j_t \neq j_{t-1} \).  At \( P_\alpha \) stage \( t \geq s \) \( P_\alpha \) acts as follows.

If \( t \) is not an eventful \( P_\alpha \) stage let \( \sigma_t \) be the partial function defined by \( \sigma_t(\pair{l_\alpha}{k_{t}}) = 0 \), that is \( \sigma_t \) asserts that \( \pair{l_\alpha}{k_{t}} \nin C \).  Enumerate into \( \axset \) (if not already present) \( \laxiom{l_\alpha +1 }{\sigma_t}{ q } \), that is place \( q \) into \( C \) if \( \pair{l_\alpha}{k_{t}} \nin C \) thereby indicating that \( P_\alpha \) completes after finite action.  If \( V_{e,t}(x) \neq \recfnl{i}{Y_1}{x} \) do nothing so that without further action we would have \( C \supfun Y_1 \).  If \( V_{e,t}(x) = \recfnl{i}{Y_1}{x} \) then enumerate \( \laxiom{l_\alpha +1 }{\sigma_t}{ y } \) into \(\axset \) as well.  This has the effect of setting \( C \supseteq Y_0 \) if \( \pair{l_\alpha}{k_{t}} \)  remains out of \(  C \).   Finally set \( k_{t+1}=k_{t} \).

If \( t \) is an eventful \( P_\alpha \) stage then enumerate the axiom \( \laxiom{l_\alpha}{\eset}{\pair{l_\alpha}{k_t}} \) into \( \axset \) and set \( k_{t+1} = k_{t} +1 \).  This has the effect of canceling the effect of any axiom enumerated by \( P_\alpha \) at any earlier stage by placing \(\pair{l_\alpha}{k_{t}} \) into \( C \).  Note that if every element of the form \( \pair{l_\alpha}{k_{t}} \) is eventually placed in \( C \) then no axiom will place \( q \) into \( C \).

Now assume that at (global) stage \( t \) \( P_\alpha \) executes it's \( s \)-th stage then define

\begin{equation*}
	f_{t}(\lh{\alpha}) = \begin{cases}
										0 & \text{ if } P_{\alpha} \text{ hasn't yet been activated }\\
										1 & \text{ if } s \text{ is an eventful stage for } P_\alpha \\
										n+2 & \text{ if } s \text{ is uneventful and there have been } n \text{ prior eventful stages.}
								 \end{cases}
\end{equation*}

If \( t \) is the first stage for which \( f_{t} \supseteq \alpha\concat[w] \) and \( t \) corresponds to \( P_\alpha \) stage \( s \) then define. 

\begin{align*}
	\begin{split}
	\delta_{\alpha\concat[w]} &= \begin{cases}
																\delta_\alpha & \text{ if } w = 0 \\
																\setcol{X\restr{r}}{\geq l_\alpha} \union \delta_\alpha \text{ where } r=\pair{l_\alpha}{k_{t}}  & \text{ otherwise }
															\end{cases}\\
	C_{\alpha\concat[w]} &= \begin{cases}
																												C_\alpha \union \delta_{\alpha\concat[w]} & \text{ if } w \neq 1 \\
																												C_\alpha \union \delta_{\alpha\concat[w]} \set{z}{z \geq k_{0}} \union \setcol{X}{l_\alpha + 1} & \text{ if } w = 1
																											\end{cases}
		\intertext{Where}
		X &= \text{ the result of } \axset_s \text{ over } C_\alpha, \delta_\alpha 
  \end{split}
\end{align*}

Note that \( \delta_{\alpha\concat[w]} \) is longer enough to restrain later modules from interfering with \( Y_0 \) and \( Y_1 \).  Also observe that for \( w > 1 \) \( \delta_{\alpha\concat[w]}  \) is defined at an uneventful \( P_\alpha \) stage \( s \) so \( \delta_{\alpha\concat[w]} \) reflects the assumption that whatever axioms are enumerated dependent on \( \pair{l_\alpha}{k_{t}} \) remain uncancelled.  On the other hand for \( w=1 \) \( \delta_{\alpha\concat[w]}  \) is defined at an eventful \( P_\alpha \) stage \( t \) so behaves as if none of the axioms dependent on \( \pair{l_\alpha}{k_{t}} \) apply.  The definition of \( C_{\alpha\concat[w]} \) guesses any unreferenced elements in column \( l_\alpha +1 \) are absent and that those in column \( l_\alpha \) are present if \( w=1 \) and absent otherwise.

\subsection{ The Basic \( \req{N}{e} \) Strategy }

We ensure that \( C \) is not computable by ensuring that \( \setcmp{\setcol{C}{l_\alpha}} \neq \setcol{\REset{e}}{l_\alpha} \).  We only place finitely many elements into \( \setcol{C}{l_\alpha} \) so if \( \setcol{\REset{e}}{l_\alpha} \) is also finite the requirement is trivially satisfied.  We ensure that if \( \setcol{\REset{e}}{l_\alpha} \) is infinite then \( \setcol{\REset{e}}{l_\alpha} \isect \setcol{C}{l_\alpha} \neq \eset \) by enumerating \( \laxiom{l_\alpha}{\eset}{y} \) into \( \axset \) at the first \( P_\alpha \) stage \( s \) and least \( y \geq s \) with  \( y \in \setcol{\omega}{l_\alpha} \) and \(  y \nin \delta_\alpha \) for which we observe \( y \in \REset[][s]{e} \).  We say \( P_\alpha \) acts at such a stage and once \( P_\alpha \) has acted we never let it do so again.

\begin{align*}
	f_{t}(\lh{\alpha}) &= \begin{cases}
										0 & \text{ if } N_\alpha \text{ hasn't acted} \\
										1 & \text{ if } N_\alpha \text{ has acted }
								 \end{cases}\\
	C_{\alpha\concat[w]} &= C_\alpha \union \setcol{X}{l_\alpha}\\
	\delta_{\alpha\concat[w]} &= \setcol{X\restr{m}}{\geq l_\alpha} \union \delta_\alpha \\
	\intertext{Where:}
	X & \text{ is the result of } \axset_s \text{ over } C_\alpha, \delta_\alpha \\
	s & \text{ is the first stage with } f_s \supseteq \alpha\concat[w] \\
	m & \text{ is larger than any number mentioned so far.}
\end{align*}

\section{Verification}

We now verify that the construction above produces the desired set \( C \).  By lemma \ref{lem:build-wrea} we have evidently built an \( \REA[\omega] \) set so all that remains is to show that \( C \nTleq \Tzero \) and \( C \Tgeq X \in \deltaZeroTwo \implies X \Tleq \Tzero \).

\begin{lemma}\label{lem:true-path-exists}
	\( f = \liminf_{s \to \infty} f_{s}  \) is well defined.  Furthermore if \( \alpha \subseteq f \) then \( P_\alpha \) is executed infinitely often.
\end{lemma}
\begin{proof}
Suppose \( n \) is the least such that the lemma fails for \( f\restr{n} = \alpha^{+} \).  Evidently \( P_{\alpha} \) can't implement \( \req{N}{e} \) as  \( f_{s}(\alpha) \) would either remain \( 0 \) or switch permanently to \( 1 \).  So assume that \( P_\alpha \) implements \( \req{R}{e,i} \).  But in this case \( f(n) \) could only be undefined if for any \( m \) there was some stage \( t \) such that if \( s > t \) \( f_{s}(n) > m \).  However, this would entail there were infinitely many eventful stages.  Hence \( \liminf_{s\to\infty } f_{s}(n) = 1 \).  The second half of the statement follows directly from the construction.
\end{proof}

\begin{lemma}\label{lem:disagree}
	Suppose \( P_\alpha \) implements \( \req{R}{e,i} \) then for all \( n \)  \( \delta_{\alpha\concat[1]} \incompat \delta_{\alpha\concat[n+2]}  \) whenever both are defined.
\end{lemma}
\begin{proof}
	Let \( q \in \setcol{\omega}{l_\alpha+1} \) be the flag location selected during the execution of \( P_\alpha \).  By the remarks at the end of subsection \ref{subsec:R_ei-strategy} \( \delta_{\alpha\concat[1]}(q)=0 \) while \( \delta_{\alpha\concat[n+2]}(q)=1  \).  
\end{proof}

\begin{lemma}\label{lem:incom-after-s}
	If \( \alpha \subseteq f, f_{s} \) then for all  \( \beta \incompat \alpha \) if \( P_\beta \) enumerates the axiom \( \laxiom{l}{\sigma}{y} \) after stage \( s \) then \( \sigma \incompat \delta_\alpha \).
\end{lemma}
\begin{proof}
	Suppose the claim holds for \( \alpha^{-} \).  If \( P_{\alpha^{-}} \) implements \( \req{N}{e} \) then for all \(  t > s \) \( f_{t} \supseteq \alpha^{-} \implies f_{t} \supseteq \alpha \) so the claim holds for \( \alpha \).  Similarly if \( P_{\alpha^{-}} \) implements \( \req{R}{e,i} \) and \( \alpha=\alpha^{-}\concat[w] \) for \( w \neq 1 \) then then claim holds for \( \alpha \).  So suppose \( \alpha=\alpha^{-}\concat[1] \).  By construction if \( {t} \geq {s} \) \( f_{t} \nsupseteq \alpha^{-}\concat[0] \).  Hence only \( \beta \) satisfying \( \beta \supseteq \alpha^{-}\concat[n+2] \) for some \( n \) are of concern.  But by the preceding lemma \( \delta_{\alpha^{-}\concat[1]} \incompat \delta_{\alpha^{-}\concat[n+2]}  \).  But if \( \beta \supseteq \alpha^{-}\concat[n+2]  \) then constraint \ref{add-const:delta_alpha} ensures that if \( P_\beta \) enumerates \( \laxiom{l}{\sigma}{y}  \) then \( \alpha^{-}\concat[n+2] \subfun \sigma \) so the lemma also holds for \( \alpha \).
\end{proof}

\begin{lemma}\label{lem:accurate}
	Suppose \( \alpha \subset f, f_s \) and \( \axset_s \) yields \( X_s \) over \( C_\beta, \delta_\beta \) with \( \beta \subset \alpha \).  Then \( X_s \supfun \delta_\alpha \) and \( \setcol{X_s}{<l_\alpha} \subset \setcol{C}{<l_\alpha} \).  Furthermore for every \( r \) there are infinitely many \( s \) such that \( X_s \supfun C\restr{r} \)
\end{lemma}
\begin{proof}
	The first claim follows by straightforward induction on \( \gamma \) with \( \beta \subseteq \gamma \subseteq \alpha \).  Since \( \setcol{X_s}{<l_\gamma} \subset \setcol{C}{<l_\gamma} \) every axiom already enumerated by \( P_\gamma \) applies in a straightforward manner as they only reference elements outside of \( \dom \delta_\gamma \) via constraint \ref{add-const:reverse}, i.e., the axioms have effect if outside of \( \delta_\gamma \) we haven't added elements not in \( C \).  By lemma \ref{lem:incom-after-s} we don't have to worry about nodes incompatible with \( \alpha \) and it is straightforward to check from the construction that when \( f_s \supseteq \gamma^{+} \) the axioms enumerated so far by \( P_\gamma \) cause \( X_s \supfun \delta_{\gamma^{+}} \).
	
	To prove the second part of the lemma simply pick \( \alpha \subset f \) so large that  \( \omega\restr{r} \subset \setcol{\omega}{< l_\alpha} \).  Now merely choose \( s \) such that \( f_s \supset \alpha \) large enough that the axioms responsible for placing every element into \( C\restr{r} \) have already been enumerated.
\end{proof}

\begin{lemma}\label{lem:satisfies-negative}
	\( C \) is not computable.
\end{lemma}
\begin{proof}
If \( C \) were computable then \( \setcmp{C} = \REset{e} \) for some \( e \).  Now pick \( \alpha \subseteq f \) such that \( P_\alpha \) implements \( \req{N}{e} \).  Now if \( \setcol{\REset{e}}{l_\alpha} \) is infinite  then there is some stage \( s \) such that \( P_\alpha \) acts to make \( \REset{e} \isect C \neq \eset \).  On the other hand if \( \setcol{\REset{e}}{l_\alpha} \) is finite then as \( \setcol{C}{l_\alpha} \) is also finite \( \setcmp{C} \neq \REset{e}  \).

\end{proof}

\begin{lemma}\label{lem:can-reverse}
	Suppose that \( \alpha \subseteq f \),  \( P_\alpha \) implements \( \req{R}{e,i} \) and there are stages \( {s_0} < {s_1} \) at which \( P_\alpha \) not yet activated with \( \alpha \subseteq f_{s_0}, f_{s_1} \) such that \( \axset_{s_j} \) yields \( C_{j} \)  over \( C_\alpha, \delta_\alpha \) for \( j=0,1 \) and \( \recfnl{i}{C_0}{}\conv[s_0] \incompat \recfnl{i}{C_1}{}\conv[s_1] \) then there are \( Y_0, Y_1 x,y \) such that \( Y_0,Y_1,x,y \) satisfy the  conditions \ref{activate-cond} at stage \( s_1 \)
\end{lemma}
\begin{proof}
	If \( Y_j = C_j\restr{s_j} \) by our convention on use we may know that \( \recfnl{i}{C_j}{}\conv[s_j]=\recfnl{i}{Y_j}{}\conv[s_j] \).  By lemma \ref{lem:incom-after-s} any axiom enumerated by \( \beta \incompat \alpha \) after \( s_0 \) will have no effect on \( C \)  and by construction the effects of all \( \beta \subsetneq \alpha \) are accounted for in \( C_\alpha, \delta_\alpha \).   As no \( \beta \supsetneq \alpha \) is allowed to affect column \( l_\alpha \) or \( l_\alpha +1 \) and \( P_\alpha \) has yet to enumerate any axioms we know that \( \setcol{Y_0}{l_\alpha} = \setcol{Y_1}{l_\alpha} = \setcol{\delta_\alpha}{l_\alpha}  \).

Now select \( y =\pair{l_\alpha +1}{m} \) where m is the first large number used by constraint \ref{add-const:reverse} for \( P_\alpha \) after stage \( s_0 \).  By constraint \ref{add-const:reverse} every axiom \( \pi \) enumerated by \( \beta \supsetneq \alpha \) after stage \( s_0 \) is enumerated dependent on sending \( y \) to \( 0 \), i.e., predicated on \( y \nin C \).  Thus, as \( y > {s_0} \) we have \( \axset_{s} \union \set{ \laxiom{l_\alpha+1}{\eset}{y} }{} \) yields some \( X \supfun Y_0 \) over \( C_\alpha, \delta_\alpha \).  The other conditions follow trivially.
\end{proof}

\begin{lemma}
	If \( V \in \deltaZeroTwo \) and \( V \Tleq C \) then \( V  \) is computable.
\end{lemma}
\begin{proof}
	Pick \( e \) such that \( V= \lim_{s\to\infty} V_{e,s} \), \( i \) such that \( \recfnl{i}{C}{} = V \) and \( \alpha \subset f \) such that \( P_\alpha \) implements \( \req{R}{e,i} \).  By construction if \( P_\alpha \) is ever activated then \( \recfnl{i}{C}{} \neq V \).  So suppose \( P_\alpha \) is never activated.  We compute \( V(x) \) as follows.  Wait for a stage \( s \) such that \( f_{s} \supseteq \alpha \) such that \( \axset_s \) yields \( Y_{s} \) over \( C_\alpha, \delta_\alpha \) and \( \recfnl{i}{Y_{s}}{x}\conv[s] \).  Use this value for \( V(x) \).

Such a stage must exist as by lemma \ref{lem:accurate} we can find \( s \) such that \( Y_{s} \) is equal to \( C \) on the use of \( \recfnl{i}{C}{x} \).  As the \( {s} \) just mentioned yields the correct value so too must our computation or there are stages \( {s_0}, {s_1} \) as in lemma \ref{lem:can-reverse} so \( P_\alpha \) is activated.  Contradiction.
\end{proof}

This completes the proof of theorem \ref{thm:no-delta2-bounding}. 

\section{Generalizations}

At this point one might naturally wonder if this result could be improved by moving to ordinals past \( \omega \).  One might conjecture there is some \( \REA[\omega \cdot \omega] \) degree \( C \) that forms a nontrivial minimal pair with \( \zerojj \).  Disappointingly this conjecture turns out to be false.  We sketch the proof below following the same approach as in lemma \ref{prop:delta3-bounding} but now considering limit stages.  The notation we use for computable ordinals is from \cite{higher-recursion-theory} and the definition of \REA[\alpha] degrees can be found in \cite{pseudojump-operators-II}.  Note that for the remainder of the paper we let \( \alpha, \beta, \lambda \) and \(  \gamma \) range over \( \kleeneO \), i.e., notations for constructive ordinals

\begin{lemma}\label{lem:indexes-for-limits}
	Suppose \( C_{\lambda} = \Tplus[\beta \kleenel \lambda][] C_{\beta} \) if \( \lambda   \) a limit, \( C_{\gamma} \kleenePlus 1 = \REset[{C_{\gamma}}][]{f(\gamma)} \Tplus C_{\gamma}  \) and \( C_{0}=\eset \).   If \( C_{\alpha} \Tleq \Tzero \) and \( f(\beta) \) is defined for all \( \beta \kleenel \alpha \) then \( \zerojj \) can (uniformly in \( \alpha \)) compute an index for \( C_\alpha \) as a \ce set.
\end{lemma}
\begin{proof}
	We prove this using definition via effective transfinite recursion.  We will define a computable function \( I(e) \) such that \( \recfnl{I(e)}{\zerojj}{\gamma} = i_\gamma \) with \( \REset{i_\gamma}=C_\gamma \) if for all \( \beta \kleenel \gamma \) \( \recfnl{e}{\zerojj}{\beta}=i_\beta \) and \( \REset{i_\beta}=C_\beta \).  Then by application of the recursion theorem \cite{on-notation-for-ordinal-numbers} we recover a fixed point \( e \) such that \( \recfnl{I(e)}{\zerojj}{}\cequiv\recfnl{e}{\zerojj}{} \) is our desired \( \zerojj \) computable function.
	
Before we construct \( I(e) \) we observe that there is a total computable function \( h \)  such that for all \( \beta \) if \( C_\beta=\REset{i}{}{} \) and \( \setcmp{C_\beta}=\REset{\hat{i}}{}{} \) then \( C_{\beta \kleenePlus 1}= \REset{h(\beta, i, \hat{i})} \).  The existence of \( h \) follows immediately from the computability of \( f \) and definition of \( C_{\beta \kleenePlus 1} \). Additionally there is a computable function \( g \) such that if \( \REset{i}=C_\gamma \) and \( \beta \kleenel \gamma \) then \( g(\gamma, \beta,i )=i' \) with \( C_{\beta}=\REset{i'}{}{} \).  As \( g \) merely unwraps some number of effective join operations it is straightforward to verify it exists.
	
	If \( \gamma=0 \) then \( \recfnl{I(e)}{\zerojj}{\gamma} \) returns a \ce index for the empty set.  If \( \gamma=\beta \kleenePlus 1 \) then \( \recfnl{I(e)}{\zerojj}{\gamma} \) first runs \( \recfnl{e}{\zerojj}{\beta} \) to extract \( i_\beta \) and then computes an index \( \hat{i_\beta} \) for the compliment of \( \REset{i_\beta} \).  The computation then returns \( h(\gamma, i_\beta, \hat{i_\beta}) \) as the index for \( C_{\gamma} \).  Finally if \( \gamma \) is a limit then \( \recfnl{I(e)}{\zerojj}{\gamma} \) searches through all pairs of indexes \( i, \hat{i} \) for complimentary \ce sets and returns the first \( i \) such that:
	
\begin{align*}
	& \forall[\beta \kleenel \gamma]\left[ \recfnl{e}{\zerojj}{\beta}=j_0 \land g(\gamma, \beta,i)=j_1 \implies \REset{j_0}=\REset{j_1}  \right] \\
	& \forall[\beta]\left[ \exists[x](\pair{\beta}{x} \in \REset{i}) \implies \beta \kleenel \gamma  \right]
\end{align*}
	
Now let \( e \) be the fixed point of \( I(e) \).  It is straightforward to trace out the definitions to verify that \( \recfnl{e}{\zerojj}{} \) behaves correctly at \( 0 \) and at every successor and limit stage so by transfinite induction \( \recfnl{e}{\zerojj}{\alpha} \) satisfies the lemma.

\end{proof}

\begin{proposition}\label{prop:generalized-zerojj-bounding}
	Suppose that \( C \) is of non-computable \( \REA[\alpha] \) degree then \( C \Tmeet \zerojj \neq \Tzero \). 
\end{proposition}
\begin{proof}
	By the definition of \( \REA[\alpha] \) sets \( C = C_{\alpha} \) where \( C_\alpha \) is defined as in lemma \ref{lem:indexes-for-limits} relative to some computable function \( f \).  Thus there is some least \( \beta \leq_{\mathcal{O}} \alpha \) such that \( C_{\beta}  \) isn't computable.  If \( \beta \) is a successor then just as in proposition \ref{prop:delta3-bounding} \( C_\beta \Tleq \zeroj  \).  So assume \( \beta \) is a limit.  By lemma \ref{lem:indexes-for-limits} we can uniformly find a \ce code for each \( C_{\gamma} \) with \( \gamma \kleenel \beta \).  To determine if \( \pair{\gamma}{x} \in C_\beta \) we first ask \( \zerojj \) if \( \gamma \kleenel \beta \).  If not \( \pair{\gamma}{x} \nin C_\beta \).  Otherwise ask \( \zerojj \) for a \ce code \( i \) for \( C_\gamma \) and report \( \pair{\gamma}{x} \in C_\beta \) iff \( \zerojj \) determines \( x \in C_\gamma \).  Hence in either case \( \Tzero \Tless C_\beta \Tleq \zerojj \) and as \( C_\beta \Tleq C  \) we have \( C \Tmeet \zerojj \neq \Tzero \). 
\end{proof}

\bibliographystyle{amsplain}
\bibliography{logic}
\end{document}